\documentstyle[amscd,amssymb,verbatim,12pt]{amsart}
\newcommand{\bfv}{\mbox{\bf v}}

\newcommand{\Id}{\mbox{\rm Id}}

\newcommand{\Ker}{\mbox{\rm Ker}}

\newcommand{\Tr}{\mbox{\rm Tr}}

\newcommand{\HH}{\mbox{\rm H}}

\newcommand{\vol}{\mbox{\rm vol}}

\newcommand{\R}{{\Bbb R}}
\newcommand{\C}{{\Bbb C}}

\newcommand{\Z}{{\Bbb Z}}

\theoremstyle{plain}
\newtheorem{definition}{Definition}

\newtheorem{lemma}{Lemma}

\newtheorem{proposition}{Proposition}

\numberwithin{equation}{section}
\renewcommand{\rm}{\normalshape}
\begin{document}
\title{Remark About Heat Diffusion on Periodic Spaces}
\author{John Lott}
\address{Department of Mathematics\\
University of Michigan\\
Ann Arbor, MI  48109-1109\\
USA}
\email{lott@@math.lsa.umich.edu}
\thanks{Research supported by NSF grant DMS-9403652}
\date{July 17, 1997}
\maketitle
\begin{abstract}
Let $M$ be a complete Riemannian manifold with a free cocompact
$\Z^k$-action.  Let $k(t, m_1, m_2)$ be the heat kernel on $M$. We compute
the asymptotics of $k(t, m_1, m_2)$ in the limit in which 
$t \rightarrow \infty$ and $d(m_1, m_2) \sim \sqrt{t}$. We show that in this
limit, the
heat diffusion is governed by an effective Euclidean metric on $\R^k$
coming from the Hodge inner product on $\HH^1(M/\Z^k; \R)$.
\end{abstract}
\section{Introduction}
Let $M$ be a complete connected oriented $n$-dimensional Riemannian manifold. 
Let $k(t, m_1, m_2)$ be the time-$t$ heat kernel on $M$.
The usual ansatz to approximate $k(t, m_1, m_2)$ is to say that
\begin{equation} \label{ansatz}
k(t, m_1, m_2) \sim P(t, m_1, m_2) \: 
e^{- \: \frac{d(m_1, m_2)^2}{4t}}
\end{equation}
where $e^{- \: \frac{d(m_1, m_2)^2}{4t}}$ is considered to be the leading
term and $P(t, m_1, m_2)$ is a correction term which can be computed
iteratively. There are results which make this
precise. For example
\cite{Azencott (1981)}, if $m_1$ and $m_2$ are nonconjugate then
as $t \rightarrow 0$,
\begin{equation}
k(t, m_1, m_2) = \sum_{\gamma}
\frac{(\det \: d(\exp_{m_1})_{v_\gamma})^{-1/2}}{(4\pi t)^{n/2}}
\: e^{- \: \frac{d(m_1, m_2)^2}{4t}} \: (1+O(t)).
\end{equation}
Here the sum is over minimal geodesics $\gamma : [0,1] \rightarrow M$ 
joining $m_1$ to $m_2$ of the form
$\gamma(s) = \exp_{m_1} (s v_\gamma)$. 
For another example, if $M$ has bounded geometry then
lower and upper
heat kernel bounds \cite{Cheeger-Yau (1981),Davies-Pang (1989)}
imply that (\ref{ansatz}) is a
good approximation if $d(m_1, m_2) > > t$, in the sense that
$- \ln(k(t, m_1, m_2))$ is well-approximated by $\frac{d(m_1, m_2)^2}{4t}$.

One can ask if the ansatz (\ref{ansatz}) is relevant for other asymptotic
regimes. In this paper we look at the case when $M$ has a periodic metric,
meaning that $\Z^k$ acts freely by orientation-preserving isometries 
on $M$, with $X = M/\Z^k$ compact.
We consider the asymptotic regime
in which $t \rightarrow \infty$ and $d(m_1, m_2) \sim \sqrt{t}$.
As the typical time-$t$ 
Brownian path will travel a distance comparable to $\sqrt{t}$, this
is the regime which contains the bulk of the diffusing heat. We show that
in this regime, (\ref{ansatz}) is no longer a valid approximation.
Instead, the heat diffusion
is governed by an effective Euclidean metric on $\R^k$. This metric
is constructed using the Hodge inner product on $\HH^1(X; \R)$.

To state the precise result,
let ${\cal F}$ be a fundamental domain in $M$ for the $\Z^k$-action.
Given $\bfv \in \Z^k$, put
\begin{equation}
k(t, \bfv) = \int_{\cal F} k(t, m, \bfv \cdot m)
\: d\vol(m). 
\end{equation}
This is independent of the choice of fundamental domain ${\cal F}$.

The covering $M \rightarrow X$ is classified by a map
$\nu : X \rightarrow B\Z^k$, defined up to homotopy, which is
$\pi_1$-surjective. It induces
a surjection $\nu_* : \HH_1(X; \R) \rightarrow \R^k$ and an injection
$\nu^* : (\R^k)^* \rightarrow \HH^1(X; \R)$. Let 
$\langle \cdot, \cdot \rangle_{H^1(X; \R)}$ be the Hodge inner product on
$\HH^1(X; \R)$. 

\begin{definition}
The inner product $\langle \cdot, \cdot \rangle_{(\R^k)^*}$ on $(\R^k)^*$ is
given by 
\begin{equation}
\langle \cdot, \cdot \rangle_{(\R^k)^*} = 
\frac{(\nu^*)^* \langle \cdot, \cdot \rangle_{H^1(X; \R)}}{\vol(X)}.
\end{equation}
The inner product $\langle \cdot, \cdot \rangle_{\R^k}$ is the dual inner
product on $\R^k$.
\end{definition}

Let $\vol(\R^k/\Z^k)$ be the volume of a lattice cell in $\R^k$, measured
with $\langle \cdot, \cdot \rangle_{\R^k}$.

\begin{proposition} \label{mainprop}
Fix $C > 0$. Then in the region 
$\{ (t, \bfv) \in \R^+ \times \Z^k : \langle \bfv, \bfv \rangle_{\R^k} 
\le C t\}$, as $t \rightarrow \infty$ we have
\begin{equation}
k(t, \bfv) = \frac{\vol(\R^k/\Z^k)}{(4 \pi t)^{k/2}} \: e^{- \: 
\langle \bfv, \bfv \rangle_{\R^k}/(4t) } \: + \: O(t^{- \: \frac{k+1}{2}})
\end{equation}
uniformly in $\bfv$.
\end{proposition}
\noindent
{\bf Examples :} \\
1. If $M = \R^k$ with a flat metric $\langle \cdot, \cdot \rangle_{flat}$
then one can check that $\langle \cdot, \cdot \rangle_{\R^k} =
\langle \cdot, \cdot \rangle_{flat}$, so one recovers the standard
heat kernel.\\
2. If $n = 2$ then $\langle \cdot, \cdot \rangle_{H^1(X; \R)}$ is 
conformally-invariant.  Hence in this case, the heat kernel asymptotics only 
depend on $\vol(X)$ and the induced complex structure on $X$.\\

One can get similar pointwise estimates on $k(t, m_1, m_2)$ by the
same methods. We omit the details.

The result of Proposition \ref{mainprop} is an example of the
phenomenon of ``homogenization'', which has been much-studied for
differential operators on $\R^n$.  Homogenization means that
in an appropriate scaling limit, the solution to a problem is governed
by the solution to a spatially homogeneous problem; see
\cite{Batty-Bratteli-Jorgesen-Robinson (1995)} and references
therein.  Thus it is not surprising that the answer in Proposition
\ref{mainprop} has a homogeneous form.  The point of the present paper is
to show how one can compute the exact asymptotics in the general
geometric setting.

We remark that when $t \rightarrow \infty$ and $d(m_1, m_2) > > t$, the 
asymptotic expression
(\ref{ansatz}) also shows homogenization. This follows from the result of
D. Burago
\cite{Burago (1992)} 
that there is a Banach norm $\parallel \cdot \parallel$ on 
$\R^k$ and a constant $c > 0$ such that if $m \in M$ and
$\bfv \in \Z^k$ then 
$\left| d(m, \bfv \cdot m) - \parallel v \parallel \right| \le c$.
Thus as $t \rightarrow \infty$, if $d(m_1, m_2) \sim \sqrt{t}$ 
then the effective
geometry is $(\R^k, \langle \cdot, \cdot \rangle_{\R^k})$, while if
$d(m_1, m_2) > > t$ then the effective geometry is $(\R^k, \parallel \cdot
\parallel)$.

It would be interesting if one could extend the results of this paper
to the setting in which $\Gamma$ is a
nonabelian discrete group, such as the fundamental group of
a closed hyperbolic surface.  In this case, the relevant scaling regime 
should be
$t \rightarrow \infty$ and $d(m_1, m_2) \sim t$, as the typical Brownian
path on the Poincar\'e disk has a constant radial velocity.

I thank the IHES for its hospitality while this work was done and
Palle Jorgensen for sending his reprints.

\section{Proof of Proposition \ref{mainprop}}

We first recall some basic facts about the eigenvalues of a
parametrized family of operators 
\cite[Chapter XII]{Reed-Simon (1978)}.

Let $M_d(\C)$ be the vector space of $d \times d$ complex matrices and 
let $M_d^{sa}(\C)$ be the subspace of self-adjoint matrices.  
Let $f : \R^k \rightarrow M_d(\C)$ be a real-analytic map.
The eigenvalues $\{\lambda_i(x)\}_{i=1}^d$ of 
$f(x)$ 
are algebraic functions of $x$, meaning the roots of a polynomial whose
coefficients are real-analytic functions of $x$, as they are
given by $\det(f(x) - \lambda) = 0$. If $\lambda_1(0)$ is a nondegenerate
eigenvalue of $f(0)$ then it extends near $x = 0$ to a real-analytic function
$\lambda_1(x)$.

If $k = 1$ and $f$ takes values in $M_d^{sa}(\C)$ then the eigenvalues
of $f$ form $d$ real-analytic functions $\{\lambda_i(x)\}_{i=1}^d$ on 
$\R$. Of course, these functions may cross.  If $k >1$ and 
$f$ takes values in $M_d^{sa}(\C)$ then it may not be true that the
eigenvalues form real-analytic functions on $\R^k$. This can be seen
in the example $f(x_1, x_2) = 
\begin{pmatrix}
0 & x_1 - i x_2 \\
x_1 + i x_2 & 0
\end{pmatrix}$. Its eigenvalues are $\pm \sqrt{x_1^2 + x_2^2}$,
which are not the union of two smooth functions on $\R^2$.
However, if $\gamma(s)$ is a real-analytic curve in $\R^2$ then
the eigenvalues of $f(\gamma(s))$ do form real-analytic functions in $s$.

If $f$ is instead an appropriate real-analytic family of operators on a
Hilbert space then one has similar results.  We refer to 
\cite[Chapter XII.2]{Reed-Simon (1978)} for the precise requirements.

To prove Proposition \ref{mainprop}, we use the method of
\cite[Section VI]{Lott (1992b)}.  The Pontryagin dual of $\Z^k$ is
$T^k = (\R^k)^*/2 \pi (\Z^k)^*$. 
Given $\theta \in T^k$, let $\rho(\theta) : \Z^k \rightarrow U(1)$
be the corresponding representation and let $E(\theta)$ be the flat
line bundle on $X$ associated to the representation 
$\pi_1(X) \stackrel{\nu_*}{\rightarrow} \Z^k 
\stackrel{\rho(\theta)}{\rightarrow} U(1)$. Let $\triangle_{\theta}$ be the
Laplacian on $L^2(X; E(\theta))$. Then Fourier analysis gives
\begin{equation} \label{eq1}
k(t, \bfv) = \int_{T^k} e^{i \theta \cdot \bfv} \: \Tr \left(e^{-t 
\triangle(\theta)} \right) \: \frac{d^k\theta}{(2\pi)^k}.
\end{equation}
Now $\Ker(\triangle(\theta)) = 0$ if $\theta \ne 0$ and
$\Ker(\triangle(0)) \cong \C$ consists of the constant functions on $X$.

In order to write all of the operators $\triangle(\theta)$ as acting on
the same Hilbert space, let $\{\tau^j\}_{j=1}^k$ be a set of harmonic
$1$-forms on $X$ which gives an integral basis of $(\Z^k)^* \subset (\R^k)^*
\subseteq \HH^1(X; \R)$.  Let $e(\tau^j)$ denote exterior multiplication by
$\tau^j$ on $C^\infty(X)$ and let $i(\tau^j)$ denote interior multiplication
by $\tau^j$ on $\Omega^1(X)$.
Putting
\begin{equation} \label{d}
d(\theta) = d + i \sum_{j=1}^k \theta_j e(\tau^j)
\end{equation}
and
\begin{equation} \label{d^*}
d^*(\theta) = d - i \sum_{j=1}^k \theta_j i(\tau^j),
\end{equation}
$\triangle(\theta)$ is
unitarily equivalent to the self-adjoint operator $d(\theta)^* d(\theta)$
(which we shall also denote by $\triangle(\theta)$)
acting on $L^2(X)$. Because $\triangle(\theta)$ is quadratic in
$\theta$, it is easy to see that $\{\triangle(\theta)\}_{\theta \in T^k}$
is an analytic family of type (A) in the sense of 
\cite[Chapter XII.2]{Reed-Simon (1978)}, so we can apply analytic eigenvalue
perturbation theory.  In particular, if $\{\lambda_i (\theta)\}_{i \in \Z^+}$ 
are the eigenvalues of $\triangle(\theta)$, arranged in increasing order and
repeated if there is a multipicity greater than one, then
$\lambda_1 (\theta) \ge 0$ and $\lambda_1 (\theta) = 0$ if and only if
$\theta = 0$, in which case it is a nondegenerate eigenvalue. Thus
$\lambda_1$ extends to a real-analytic function in a neighborhood of 
$\theta = 0$. So for sufficiently small $\epsilon > 0$, there is a neighborhood
$U \subseteq T^k$ of $0 \in T^k$ such that \\
1. If $\theta \notin U$ then
$\lambda_1 (\theta) > \epsilon$. \\
2. Restricted to $U$, $\lambda_1$ is a 
real-analytic function which represents a nondegenerate
eigenvalue and $\lambda_2 > \epsilon$. \\

From (\ref{eq1}), we have
\begin{equation}
k(t, \bfv) = \int_{T^k} e^{i \theta \cdot \bfv} \: 
\sum_{i=1}^\infty e^{-t \lambda_i (\theta)} \: \frac{d^k\theta}{(2\pi)^k}.
\end{equation}
Then it is easy to show that
\begin{equation}
k(t, \bfv) = \int_{U} e^{i \theta \cdot \bfv} \: 
e^{-t \lambda_1 (\theta)} \: \frac{d^k\theta}{(2\pi)^k} + O(e^{-\epsilon
t/2}),
\end{equation}
uniformly in $\bfv$.

\begin{lemma} \label{lemma1}
The Taylor's series of $\lambda_1 (\theta)$ near $\theta = 0$ starts off
as 
\begin{equation}
\lambda_1 (\theta) = \langle \theta, \theta \rangle_{_{(\R^k)^*}} + 
O(|\theta|^3).
\end{equation}
\end{lemma}
\begin{pf}
It suffices to compute $\frac{d\lambda_1(s \vec{w})}{ds} \big|_{s=0}$ and
$\frac{d^2\lambda_1(s \vec{w})}{ds^2} \big|_{s=0}$ for
all $\vec{w} \in (\R^k)^*$.
For simplicity, denote $\triangle(s \vec{w})$ by
$\triangle(s)$ and $\lambda_1(s \vec{w})$ by $\lambda(s)$.
As $\lambda(s)$ is nonnegative and $\lambda(0) = 0$, we must
have $\lambda^\prime(0) = 0$.
Let $\psi(s)$ denote a nonzero eigenfunction with eigenvalue $\lambda(s)$;
we can assume that it is real-analytic in $s$ with
$\psi(0) = 1$. 
Differentiation of
$\triangle(s) \psi(s) = \lambda(s) \psi(s)$ gives
\begin{equation} \label{eq2}
\triangle^\prime(0) \psi(0) + \triangle(0) \psi^\prime (0) = 0
\end{equation}
and
\begin{equation} \label{eq3}
\triangle^{\prime \prime}(0) \psi(0) + 2 \triangle^\prime(0) \psi^\prime(0)
+ \triangle(0) \psi^{\prime \prime}(0) = 
\lambda^{\prime \prime}(0) \psi(0).
\end{equation}
Taking the inner product of (\ref{eq3}) with $\psi(0)$ gives
\begin{equation} \label{eq4}
\langle \psi(0), \triangle^{\prime \prime}(0) \psi(0) \rangle
+ 2 \langle \psi(0), \triangle^\prime(0) \psi^\prime(0) \rangle
= 
\lambda^{\prime \prime}(0) \langle \psi(0), \psi(0) \rangle.
\end{equation}
Let $G$ be the Green's operator for $\triangle(0)$. From (\ref{eq2}),
\begin{equation} \label{eq5}
\psi^\prime(0) = c \psi(0) - G \triangle^\prime(0) \psi(0)
\end{equation}
for some constant $c$. Changing $\psi(s)$ to $e^{-cs} \psi(s)$, we may assume
that $c = 0$. 
Substituting (\ref{eq5}) into (\ref{eq4}) gives
\begin{equation} \label{eq6}
\langle \psi(0), \triangle^{\prime \prime}(0) \psi(0) \rangle
- 2 \langle \psi(0), \triangle^\prime(0) G \triangle^\prime(0) \psi(0) \rangle
= \lambda^{\prime \prime}(0) \langle \psi(0), \psi(0) \rangle.
\end{equation}

It remains to compute 
$\langle \psi(0), \triangle^{\prime \prime}(0) \psi(0) \rangle$ and
$\langle \psi(0), \triangle^\prime(0) G \triangle^\prime(0) \psi(0) \rangle$.
Put $D(s) = d_{s \vec{w}}$ and $D^*(s) =  d_{s \vec{w}}^*$.
Then $\triangle(s) = D^*(s) D(s)$.
From (\ref{d}) and (\ref{d^*}), 
$D(s)$ and $D^*(s)$ are linear in $s$, with
\begin{equation}
D^\prime(0) = i \sum_{j=1}^k w_j \:  e(\tau^j)
\end{equation}
and
\begin{equation}
(D^*)^\prime(0) = - i \sum_{j=1}^k w_j \: i(\tau^j).
\end{equation}
Then 
\begin{align} \label{eq7}
\langle \psi(0), \triangle^{\prime \prime}(0) \psi(0) \rangle & =
2 \langle \psi(0), (D^*)^\prime(0) D^\prime(0) \psi(0) \rangle \\
& = 2 \big| D^\prime(0) \psi(0) \big|^2_{H^1(X; \C)} \notag \\
& = 2 \big| \sum_{j=1}^k w_j \: \tau^j \big|^2_{H^1(X; \C)}. \notag
\end{align}
Now
\begin{align} \label{eq8}
\triangle^\prime(0) \psi(0) & = \left[ (D^*)^\prime(0) D(0) +
D^*(0) D^\prime(0) \right] \psi(0) \\
& = d^* \left( - i \sum_{j=1}^k w_j \: \tau^j \right) = 0.
\notag
\end{align}
Substituting (\ref{eq7}) and (\ref{eq8}) into (\ref{eq6}) and using the
fact that $\langle \psi(0), \psi(0) \rangle = \vol(X)$, the lemma follows.  
\end{pf}

Continuing with the proof of Proposition \ref{mainprop}, by Morse theory
and Lemma \ref{lemma1},
we can find a change of coordinates near $0 \in T^k$ with respect to which
$\lambda_1$ becomes quadratic. That is, if $B_r(0)$ denotes the ball of
radius $r$ in $(\R^k)^*$, 
we can find an $r > 0$, a neighborhood $U$ of $0 \in T^k$ and a 
diffeomorphism $\phi : B_r(0) \rightarrow U$ such that $\phi(0) = 0$,
$d \phi_0 = \Id$ and 
$\lambda_1 (\phi(x)) = \langle x, x \rangle_{(\R^k)^*}$. Then there is
some $\alpha > 0$ such that as $t \rightarrow \infty$,
\begin{equation}
k(t, \bfv) = \int_{B_r(0)} e^{i \phi(x) \cdot \bfv} \: 
e^{-t \langle x, x \rangle_{(\R^k)^*}} \: 
\det (d\phi_x) \:
\frac{d^kx}{(2\pi)^k} + O(e^{- \alpha t}),
\end{equation}
uniformly in $\bfv$.
Multiplying by a cutoff function on $(\R^k)^*$, we can write
\begin{align} \label{eqnstart}
k(t, \bfv) & = \int_{(\R^k)^*} e^{i \phi(x) \cdot \bfv} \: 
e^{-t \langle x, x \rangle_{(\R^k)^*}} \: 
g(x) \: 
\frac{d^kx}{(2\pi)^k} + O(e^{- \alpha^\prime t}) \\
& = t^{- \: \frac{k}{2}}
\int_{(\R^k)^*} e^{i \phi(\frac{x}{\sqrt{t}}) \cdot \bfv} \: 
e^{- \langle x, x \rangle_{(\R^k)^*}} \: 
g \left(\frac{x}{\sqrt{t}} \right) \: 
\frac{d^kx}{(2\pi)^k} + O(e^{- \alpha^\prime t}) \notag
\end{align}
for some $g \in C^\infty_0 \left( (\R^k)^* \right)$ 
with $g(0) = 1$ and some $\alpha^\prime > 0$. 
(Here $\phi$ has been extended to become a map
$\phi : (\R^k)^* \rightarrow (\R^k)^*$ which is the identity outside
of a compact set.)

We have now reduced to a stationary-phase-type integral.  Let
\begin{equation}
g(x) = 1 + (\nabla g) (0) \cdot x + E(x)
\end{equation}
be the beginning of the Taylor's expansion of $g$.
We can write
\begin{align}
& t^{- \: \frac{k}{2}}
\int_{(\R^k)^*} e^{i \phi(\frac{x}{\sqrt{t}}) \cdot \bfv} \: 
e^{- \langle x, x \rangle_{(\R^k)^*}} \: 
g \left(\frac{x}{\sqrt{t}} \right) \: 
\frac{d^kx}{(2\pi)^k} = \\
& \hspace{1.0in} t^{- \: \frac{k}{2}}
\int_{(\R^k)^*} e^{i \frac{x}{\sqrt{t}} \cdot \bfv} \: 
e^{- \langle x, x \rangle_{(\R^k)^*}} \: 
\left[ 1 + (\nabla g) (0) \cdot \frac{x}{\sqrt{t}} + 
E \left(\frac{x}{\sqrt{t}} \right) \right] \: 
\frac{d^kx}{(2\pi)^k} + \notag \\
& \hspace{1.0in} t^{- \: \frac{k}{2}}
\int_{(\R^k)^*} 
e^{i\frac{x}{\sqrt{t}} \cdot \bfv}
\left[ e^{i \left( \phi(\frac{x}{\sqrt{t}}) - \frac{x}{\sqrt{t}} \right) 
\cdot \bfv} - 1 \right] \: 
e^{- \langle x, x \rangle_{(\R^k)^*}} \: 
g \left(\frac{x}{\sqrt{t}} \right) \: 
\frac{d^kx}{(2\pi)^k}. \notag
\end{align}

Recall that the measure $\frac{d^kx}{(2\pi)^k}$ on $(\R^k)^*$ derives
from the product measure on $T^k = \left( \R^*/2 \pi \Z^* \right)^k$. 
Let $\langle \cdot, \cdot \rangle_{prod}$ be the standard product Euclidean
metric on $\left( \R^* \right)^k$. Let $Q$ be the self-adjoint operator
on $(\R^k)^*$
such that $\langle x, x \rangle_{(\R^k)^*} = \langle x, Qx \rangle_{prod}$.
Then a standard calculation gives
\begin{equation}
t^{- \: \frac{k}{2}}
\int_{(\R^k)^*} e^{i \frac{x}{\sqrt{t}} \cdot \bfv} \: 
e^{- \langle x, x \rangle_{(\R^k)^*}} \: 
\frac{d^kx}{(2\pi)^k} = \frac{(\det \: Q)^{-1/2}}{(4 \pi t)^{k/2}} \: e^{- \: 
\langle \bfv, \bfv \rangle_{\R^k}/(4t) }.
\end{equation}
On the other hand,
\begin{equation}
(\det \: Q)^{-1/2} = \vol(\R^k/\Z^k).
\end{equation}

By symmetry,
\begin{equation}
t^{- \: \frac{k}{2}}
\int_{(\R^k)^*} e^{i \frac{x}{\sqrt{t}} \cdot \bfv} \: 
e^{- \langle x, x \rangle_{(\R^k)^*}} \: 
(\nabla g) (0) \cdot \frac{x}{\sqrt{t}} \: 
\frac{d^kx}{(2\pi)^k}
\end{equation}
Let $c > 0$ be such that $|E(x)| \le c
\langle x, x \rangle_{(\R^k)^*}$ for all $x \in (\R^k)^*$. 
Then
\begin{equation}
\left| \int_{(\R^k)^*} e^{i \frac{x}{\sqrt{t}} \cdot \bfv} \: 
e^{- \langle x, x \rangle_{(\R^k)^*}} \: 
E \left(\frac{x}{\sqrt{t}} \right) \: 
\frac{d^kx}{(2\pi)^k} \right| \le \frac{c}{t} \:
\int_{(\R^k)^*} \langle x, x \rangle_{(\R^k)^*} \:
e^{- \langle x, x \rangle_{(\R^k)^*}} \: 
\frac{d^kx}{(2\pi)^k}.
\end{equation}

Finally,
\begin{align}
& \left| \int_{(\R^k)^*} 
e^{i\frac{x}{\sqrt{t}} \cdot \bfv}
\left[ e^{i \left( \phi(\frac{x}{\sqrt{t}}) - \frac{x}{\sqrt{t}} \right) 
\cdot \bfv} - 1 \right] \: 
e^{- \langle x, x \rangle_{(\R^k)^*}} \: 
g \left(\frac{x}{\sqrt{t}} \right) \: 
\frac{d^kx}{(2\pi)^k} \right| \le \\
& \hspace{1.0in} \parallel g \parallel_\infty \int_{(\R^k)^*} 
2 \: \left| \sin \left( \frac12 \left[ \phi \left( \frac{x}{\sqrt{t}} \right) 
- \frac{x}{\sqrt{t}} \right] \cdot \bfv
\right) \right| \:
e^{- \langle x, x \rangle_{(\R^k)^*}} \:
\frac{d^kx}{(2\pi)^k}. \notag
\end{align}
We can find a constant $c^\prime > 0$ such that 
\begin{equation}
2 \left| \sin \left( \frac12 \left[\phi (x) 
- x \right] \cdot \bfv \right) \right| \le c^\prime \:
\langle x, x \rangle_{(\R^k)^*} \: \parallel \bfv \parallel_{\R^k}
\end{equation} 
for all $x \in (\R^k)^*$ and $\bfv \in \Z^k$. 
Then
\begin{align}
&\parallel g \parallel_\infty \int_{(\R^k)^*} 
2 \: \left| \sin \left( \frac12 \left[ \phi \left( \frac{x}{\sqrt{t}} \right) 
- \frac{x}{\sqrt{t}} \right] \cdot \bfv
\right) \right| \:
e^{- \langle x, x \rangle_{(\R^k)^*}} \:
\frac{d^kx}{(2\pi)^k}
 \le \\
& \hspace{1.0in} \frac{c^\prime}{\sqrt{t}} \: 
\frac{\parallel \bfv \parallel_{\R^k}}{\sqrt{t}}
\: \parallel g \parallel_\infty \int_{(\R^k)^*} 
\langle x, x \rangle_{(\R^k)^*} \:
e^{- \langle x, x \rangle_{(\R^k)^*}} \:
\frac{d^kx}{(2\pi)^k}. \notag
\end{align}
By assumption, 
\begin{equation} \label{eqnend}
\frac{\parallel \bfv \parallel_{\R^k}}{\sqrt{t}} \le \sqrt{C}.
\end{equation}
The proposition follows from combining equations
(\ref{eqnstart})-(\ref{eqnend}).

\end{document}